# The instability of naked singularities in the gravitational collapse of a scalar field

By Demetrios Christodoulou

## 1. Introduction

One of the fundamental unanswered questions in the general theory of relativity is whether "naked" singularities, that is singular events which are visible from infinity, may form with positive probability in the process of gravitational collapse. The conjecture that the answer to this question is in the negative has been called "cosmic censorship." The present paper, which is a continuation of the work in [1], [2], [3] and [4] addresses this question in the context of the spherical gravitational collapse of a scalar field.

The problem of a spherically symmetric self-gravitating scalar field is formulated in terms of a 2-dimensional quotient space-time manifold $\mathcal{Q}$ with boundary (see [1]). The boundary of $\mathcal{Q}$ corresponds to the set of fixed points of the group action, the center of symmetry, which is a timelike geodesic $\Gamma$. The manifold $\mathcal{Q}$ is endowed with a Lorentzian metric $g_{ab}$, an area radius function $r$ and a wave function $\phi$ satisfying the following nonlinear system of partial differential equations:

$$(1.1a) \qquad r\nabla_a\nabla_b r = \frac{1}{2}g_{ab}(1 - \partial^c r\partial_c r) - r^2 T_{ab}$$

$$T_{ab} = \partial_a\phi\partial_b\phi - \frac{1}{2}g_{ab}\partial^c\phi\partial_c\phi$$

$$(1.1b) \qquad \nabla^a(r^2\partial_a\phi) = 0.$$

These imply the following equation for the Gauss curvature of $\mathcal{Q}$:

$$(1.1c) \qquad K = r^{-2}(1 - \partial^a r\partial_a r) + \partial^a\phi\partial_a\phi.$$

The mass function $m$ is defined by:

$$(1.2) \qquad 1 - \frac{2m}{r} = g^{ab}\partial_a r\partial_b r.$$

In [2] it was shown that given an initial future light cone with vertex at the center of symmetry and with a region bounded by two spheres such that the ratio of the mass contained in the region to the largest radius is large in comparison to the ratio of the radii minus 1, then a trapped region, namely a region where the future light cones have negative expansion, forms in the future terminating at a strictly spacelike singular boundary. The trapped region contains a sphere whose mass is bounded from below by a positive number depending only on the two initial radii. The results of [2] shall be used in an essential way in the present paper.

Solutions with initial data of bounded variation were considered in [1] and a sharp sufficient condition on the initial data was found for the avoidance of singularities, namely that the total variation be sufficiently small, as well as another condition implying the formation of singularities, complementing the results in [2]. Also, a sharp extension criterion was established, namely that if the ratio of mass to radius of spheres tends to zero as we approach a point at the center of symmetry from its causal past, then the solution extends as a regular solution to include a full neighborhood of the point. The structure of bounded variation solutions was studied and it was shown that at each point in the center of symmetry the solutions are locally scale invariant. Also, the behaviour of the solutions at the singular boundary was analyzed. The present paper relies on the results of [1] as well.

In [3] it was shown that when the final Bondi mass, that is, the infimum of the mass at future null infinity, is different from zero, a black hole forms of mass equal to the final Bondi mass surrounded by vacuum. The rate of growth of the redshift of light seen by faraway observers was determined and the asymptotic wave behaviour at timelike infinity and along the event horizon, the boundary of the past of future null infinity, was analyzed.

In [4] we constructed examples of solutions corresponding to regular asymptotically flat initial data which develop singularities which are not preceeded by a trapped region but have future light cones expanding to infinity. Thus naked singularities do, in fact, occur in the spherical gravitational collapse of a scalar field.

The present paper nevertheless supports the cosmic censorship conjecture. For, we shall show in the following that in the space of initial conditions the subset of initial conditions leading to the formation of naked singularities has, in a certain sense, positive codimension, consequently the occurence of naked singularities is an unstable phenomenon in the context of the spherical self-gravitating scalar field model.

We shall be concerned here with the part of the space-time manifold which lies in the past of the apparent horizon, the past boundary of the trapped region. This part includes the causal past of $\Gamma$ (see [1]). Consequently, in the region of interest the gradient of the function $r$ is a spacelike outward- directed vectorfield and we may use $r$ as a coordinate. We also use a null coordinate $u$ which is constant along the future-directed null curve from each point on $\Gamma = \partial \mathcal{Q}$ and increases toward the future. In terms of the *Bondi coordinates u* and $r$ the metric assumes the form:

$$(1.3) \qquad g_{ab}dx^a dx^b = -e^{2\nu}du^2 - 2e^{\nu+\lambda}dudr.$$

We shall express the system of equations (1.1a), (1.1b) in Bondi coordinates. It is advantageous to use the pair of null future-directed vectorfields $n$ and $l$,

$$n = 2e^{-\nu}\frac{\partial}{\partial u} - e^{-\lambda}\frac{\partial}{\partial r} \quad , \qquad l = e^{-\lambda}\frac{\partial}{\partial r}$$
$$(g(n,l) \quad = \quad -2),$$

the integral curves of which are the outgoing and incoming null curves, respectively. We have:

$$lr = e^{-\lambda}, \quad nr = -e^{-\lambda}.$$

We note that the metric function $\lambda$ has an invariant geometric meaning, for

$$\frac{4}{r^2}e^{-2\lambda}$$

is the square of the length of the mean curvature vector of the spheres which are the orbits of the rotation group in the 4-dimensional space-time manifold. It is also advantageous to express the derivatives of the wave function in terms of

$$(1.4) \qquad \theta = r\left(\frac{l\phi}{lr}\right), \quad \zeta = r\left(\frac{n\phi}{nr}\right).$$

We have

$$(1.5) \qquad \theta = r\frac{\partial\phi}{\partial r}, \quad \zeta = -2re^{\lambda-\nu}\frac{\partial\phi}{\partial u} + r\frac{\partial\phi}{\partial r}.$$

The mass function $m$, defined, in general, by (1.2), is, in terms of Bondi coordinates, given by:

$$(1.6) \qquad 1 - \frac{2m}{r} = e^{-2\lambda}.$$

The components of the energy tensor $T_{ab}$ are:

$$T(n,n) = (n\phi)^2, \quad T(l,l) = (l\phi)^2$$
$$T(n,l) = -\mathrm{tr}T = 0.$$

The trace of equation (1.1a) is

$$(1.7) \qquad r\left(\frac{\partial\nu}{\partial r} - \frac{\partial\lambda}{\partial r}\right) = e^{2\lambda} - 1$$

while the trace-free part of equation (1.1a) reduces to the following pair of equations for $m$:

$$nm = -(1/2)r^2(lr)T(n,n) = -(1/2)r^2e^{-\lambda}(n\phi)^2$$

$$lm = -(1/2)r^2(nr)T(l,l) = (1/2)r^2 e^{-\lambda}(l\phi)^2;$$

that is,

(1.8a) $$2e^{\lambda-\nu}\frac{\partial m}{\partial u} - \frac{\partial m}{\partial r} = -\frac{1}{2}e^{-2\lambda}\zeta^2$$

(1.8b) $$\frac{\partial m}{\partial r} = \frac{1}{2}e^{-2\lambda}\theta^2.$$

By virtue of (1.7) this last equation can be written as:

(1.9) $$r\left(\frac{\partial \nu}{\partial r} + \frac{\partial \lambda}{\partial r}\right) = \theta^2.$$

Finally, the wave equation (1.1b) takes in Bondi coordinates the form:

$$-2\left(\frac{\partial^2 \phi}{\partial u \partial r} + \frac{1}{r}\frac{\partial \phi}{\partial u}\right) + e^{\nu-\lambda}\left[\frac{\partial^2 \phi}{\partial r^2} + \left(\frac{2}{r} + \frac{\partial \nu}{\partial r} - \frac{\partial \lambda}{\partial r}\right)\frac{\partial \phi}{\partial r}\right] = 0.$$

In view of (1.7), the wave equation is equivalent to the following pair of equations for $\theta, \zeta$:

(1.10a) $$r\left(2e^{\lambda-\nu}\frac{\partial \theta}{\partial u} - \frac{\partial \theta}{\partial r}\right) = (e^{2\lambda} - 1)\theta + \zeta$$

(1.10b) $$r\frac{\partial \zeta}{\partial r} = -(e^{2\lambda} - 1)\zeta - \theta.$$

Thus $\phi$ is eliminated and we have a first-order system in the unknowns $\lambda, \nu, \theta, \zeta$.

Given initial data of bounded variation (see [1]) on a future light cone $C_0^+$ with vertex at the center of symmetry), we consider the maximal sphere $S_0$ in $C_0^+$ such that for each sphere $S$ in $C_0^+$ within $S_0$ there is a point $P$ on the central geodesic $\Gamma$ whose past light cone $C_P^-$ intersects $C_0^+$ at $S$, and we have a development of bounded variation in the region bounded by $C_P^-$ and $C_0^+$. Let $C_0^-$ be the incoming null hypersurface through $S_0$ in the future of $C_0^+$. Then, as we have shown in [1], the solution extends to $C_0^-$, however $C_0^-$ cannot terminate at a regular vertex on $\Gamma$. We can choose the coordinate $u$ so that

$$u = -2r$$

along the incoming null curve corresponding to $C_0^-$ in the quotient space-time $\mathcal{Q}$. The incoming null curves satisfy the equation:

(1.11) $$\frac{dr}{du} = -\frac{1}{2}e^{\nu-\lambda}.$$

Therefore we have

(1.12) $$e^{\nu-\lambda}\Big|_{u=-2r} = 1.$$

The origin $(0,0)$ in the coordinates $u, r$ set up in this way cannot correspond to a regular point. Let the sphere $S_0$ correspond to $r = a$ ($u = -2a$). We define dimensionless coordinates $t, s$ by:

$$(1.13) \qquad u = -2ae^{-t}, \quad r = ae^{s-t}.$$

Then $t$ is constant along the outgoing null curves while

$$(1.14) \qquad \frac{ds}{dt} = \beta, \qquad \beta = 1 - e^{\nu - \lambda - s},$$

along the incoming null curves. Since $s = 0$ corresponds to $C_0^-$,

$$(1.15) \qquad \beta\big|_{s=0} = 0.$$

The future of $C_0^+$ corresponds to $t > 0$, and the interior of $C_0^-$ corresponds to $s < 0$, the exterior to $s > 0$. We have:

$$(1.16a) \qquad r\frac{\partial}{\partial r} = \frac{\partial}{\partial s},$$

$$(1.16b) \qquad r\left(2e^{\lambda-\nu}\frac{\partial}{\partial u} - \frac{\partial}{\partial r}\right) = \frac{1}{(1-\beta)}\left(\frac{\partial}{\partial t} + \beta\frac{\partial}{\partial s}\right),$$

while

$$(1.16c) \qquad \frac{\partial}{\partial t} = -u\frac{\partial}{\partial u} - r\frac{\partial}{\partial r}.$$

With

$$(1.17) \qquad \kappa = e^{2\lambda},$$

equations (1.6), (1.8a), 1.8b), (1.10a), (1.10b) take in the dimensionless coordinates of the form:

$$(1.18a) \qquad \frac{\partial\kappa}{\partial t} + \beta\frac{\partial\kappa}{\partial s} = (1-\beta)\kappa(\kappa - 1 - \zeta^2)$$

$$(1.18b) \qquad \frac{\partial\kappa}{\partial s} = \kappa(1 - \kappa + \theta^2)$$

$$(1.18c) \qquad \frac{\partial\beta}{\partial s} = (1-\beta)(2-\kappa)$$

$$(1.18d) \qquad \frac{\partial\theta}{\partial t} + \beta\frac{\partial\theta}{\partial s} = (1-\beta)[(\kappa-1)\theta + \zeta]$$

$$(1.18e) \qquad \frac{\partial\zeta}{\partial s} = -(\kappa-1)\zeta - \theta.$$

Let us denote by a subscript 0 the restriction to $s = 0$. Then by virtue of (1.15), that is

$$(1.19) \qquad \beta_0 = 0,$$

the restrictions of equations (1.18a) and (1.18d) to $s = 0$ are:

$$(1.20a) \qquad \frac{d\kappa_0}{dt} = \kappa_0(\kappa_0 - 1 - \zeta_0^2)$$

$$(1.20b) \qquad \frac{d\theta_0}{dt} = (\kappa_0 - 1)\theta_0 + \zeta_0.$$

Equations (1.18b) and (1.18c) are equivalent to equations (1.9) and (1.7) for $\nu$ and $\lambda$, which in the dimensionless coordinates read:

$$(1.21a) \qquad \frac{\partial\nu}{\partial s} + \frac{\partial\lambda}{\partial s} = \theta^2$$

$$(1.21b) \qquad \frac{\partial\nu}{\partial s} - \frac{\partial\lambda}{\partial s} = e^{2\lambda} - 1.$$

By (1.14), (1.17) and (1.19),

$$(1.22) \qquad \nu_0 = \lambda_0 = \frac{1}{2}\log\kappa_0.$$

Hence integrating (1.21a) from $s = 0$ yields:

$$(1.23a) \qquad (\nu + \lambda)(t, s) = \log\kappa_0(t) + \int_0^s \theta^2(t, s')ds'.$$

Also, writing (1.21b) in the form

$$\frac{\partial e^{\nu-\lambda}}{\partial s} = e^{\nu+\lambda} - e^{\nu-\lambda},$$

or

$$\frac{\partial e^{\nu-\lambda-s}}{\partial s} = e^{\nu+\lambda+s}$$

and integrating from $s = 0$ yields:

$$(1.23b) \qquad e^{(\nu-\lambda)(t,s)+s} = 1 + \int_0^s e^{(\nu+\lambda)(t,s')+s'}ds'.$$

Using (1.21b) we can write (1.18e) in the form:

$$\frac{\partial(e^{\nu-\lambda}\zeta)}{\partial s} = -e^{\nu-\lambda}\theta.$$

Hence, integrating from $s = 0$ we obtain

$$(1.24a) \qquad (e^{\nu-\lambda}\zeta)(t, s) = \zeta_0(t) + \xi(t, s)$$

where

$$(1.24b) \qquad \xi(t, s) = -\int_0^s (e^{\nu-\lambda}\theta)(t, s')ds'.$$

Let us define the mass ratio

$$(1.25) \qquad \mu = \frac{2m}{r}.$$

We then have,

$$(1.26) \qquad \kappa = \frac{1}{1-\mu}$$

(see (1.6)). The fact that $\mu$ is nonnegative (see [1]) implies:

$$(1.27a) \qquad \kappa \geq 1.$$

In particular,

$$(1.27b) \qquad \kappa_0 \geq 1.$$

## 2. The first instability theorem

In the following we confine attention to the exterior of $C_0^-$ and the future of $C_0^+$ : $s, t > 0$. Let us define

$$(2.1) \qquad \gamma(t) = \int_0^t (\kappa_0(t') - 1) dt'.$$

LEMMA 1.

$$\kappa_0(t) \leq 2\kappa_0(0) e^{\gamma(t)}.$$

*Proof.* According to (1.8a), $m$ is nonincreasing along incoming null curves. Consequently, by (1.13) and (1.25) the function $\mu_0 e^{-t}$ is nondecreasing; hence, by (1.26),

$$(2.2a) \qquad \kappa_0(t) \leq \frac{1}{1 - \left(1 - \frac{1}{\kappa_0(0)}\right) e^t}$$

provided that

$$e^t < \frac{1}{1 - \frac{1}{\kappa_0(0)}}.$$

On the other hand, since $t' \leq t$ implies

$$\mu_0(t') e^{-t'} \geq \mu_0(t) e^{-t},$$

we also have that

$$\kappa_0(t') - 1 = \frac{\mu_0(t')}{1 - \mu_0(t')} \geq \frac{\mu_0(t) e^{t'-t}}{1 - \mu_0(t) e^{t'-t}}.$$

Hence,

$$\begin{aligned} \gamma(t) &\geq \int_0^t \frac{\mu_0(t) e^{t'-t}}{1 - \mu_0(t) e^{t'-t}} dt' \\ &= \log\left[\frac{1 - \mu_0(t) e^{-t}}{1 - \mu_0(t)}\right] > \log\left[\frac{1 - e^{-t}}{1 - \mu_0(t)}\right] \end{aligned}$$

since $\mu_0 < 1$; that is,

(2.2b) $$\kappa_0(t) < \frac{e^{\gamma(t)}}{1 - e^{-t}}.$$

By setting

$$t_1 = \log\left[\frac{1 - \frac{1}{2\kappa_0(0)}}{1 - \frac{1}{\kappa_0(0)}}\right],$$

we see that (2.2a) holds for $t \in [0, t_1]$, which yields

$$\kappa_0(t) \leq \frac{1}{1 - \left(1 - \frac{1}{\kappa_0(0)}\right)e^{t_1}} = 2\kappa_0(0)$$

for $t \in [0, t_1]$, while by (2.2b)

$$\kappa_0(t) < \frac{e^{\gamma(t)}}{1 - e^{-t_1}} = (2\kappa_0(0) - 1)e^{\gamma(t)}$$

for $t \in (t_1, \infty)$.

Since $\gamma$ is a nondecreasing function, either $\gamma$ is bounded, in which case it tends to a finite limit $\gamma(\infty)$ as $t \to \infty$, or $\gamma$ is unbounded, in which case $\gamma(t) \to \infty$ as $t \to \infty$.

LEMMA 2.   *If $\gamma$ is bounded then $\mu_0(t) \to 0$ as $t \to \infty$.*

*Proof.* Let $t > t_0$. Following the proof of Lemma 1 we obtain

$$\begin{aligned}
\gamma(t) - \gamma(t_0) &\geq \int_{t_0}^{t} \frac{\mu_0(t)e^{t'-t}}{1 - \mu_0(t)e^{t'-t}}dt' \\
&= \log\left[\frac{1 - \mu_0(t)e^{t_0-t}}{1 - \mu_0(t)}\right];
\end{aligned}$$

that is,

$$\frac{\mu_0(t)(1 - e^{t_0-t})}{1 - \mu_0(t)} \leq e^{\gamma(t)-\gamma(t_0)} - 1$$

which implies

(2.3) $$\mu_0(t) \leq \frac{e^{\gamma(t)-\gamma(t_0)} - 1}{1 - e^{t_0-t}}.$$

The result follows by setting $t_0 = t - 1$ in (2.3), in view of the fact that $\gamma$ bounded implies $\gamma(t) - \gamma(t-1) \to 0$ as $t \to \infty$.

As we have noted, the incoming null hypersurface $C_0^-$ cannot terminate at a regular vertex on $\Gamma$. The past of $C_0^-$ is a *terminal indecomposable past set* in the terminology of [5]. In view of Lemma 2 and the extension criterion mentioned in the introduction, we shall assume in the following that the function $\gamma$ is unbounded. Defining another null coordinate $v$ which is constant along the incoming null curves and satisfies

$$v = 2(r - a)$$

along the outgoing null curve corresponding to $C_0^+$, we see that $\mathcal{Q}$ becomes a domain in the $u, v$ plane in which $u < 0$, $C_0^-$ corresponds to the line $v = 0$ in $\mathcal{Q}$, and the terminal indecomposable past set to the origin $O$ which lies on the boundary of $\mathcal{Q}$ in the $u, v$ plane. The point $O$ is the past end point of the *central component* $\mathcal{B}_0$ of the *singular boundary* $\mathcal{B}$ of $\mathcal{Q}$ (see [1]). The function $r$ extends continuously to $O$ where it vanishes.

The *apparent horizon* $\mathcal{A}$ is the set of points of $\mathcal{Q}$ at which $\partial r/\partial v = 0$. Each point of $\mathcal{A}$ corresponds to a sphere which has maximal area in the future light cone with vertex on $\Gamma$ in which it is contained. According to the results of [1], $\mathcal{A}$, if nonempty, is a spacelike curve which may contain outgoing null segments but does not contain incoming null segments. In fact $\mathcal{A}$ is given by

$$(2.4a) \qquad \mathcal{A} = \{(u, v_0(u)) \, : \, u \in (\underline{u}^*, 0)\} \bigcup \left( \bigcup_n \{u_n\} \times I_n \right),$$

where $v_0$ is a strictly decreasing function in $(\underline{u}^*, 0)$, $\underline{u}^* > -2a$, and the intervals

$$(2.4b) \qquad I_n = \left( \lim_{u \to u_n^+} v_0(u), \ \lim_{u \to u_n^-} v_0(u) \right)$$

correspond to the (denumerable) points of discontinuity of $v_0$. Also,

$$(2.4c) \qquad v_0(u) \to \infty \ \text{as} \ u \to \underline{u}^*.$$

The future light cone, with vertex on $\Gamma$, which corresponds to the outgoing null curve $u = \underline{u}^*$ is the *event horizon* $\mathcal{H}$ (see [2]).

Moreover, $\mathcal{A}$ can equivalently be defined as the set of points of $\mathcal{Q}$ at which $\mu = 1$. The past of $\mathcal{A}$ in $\mathcal{Q}$, the domain of the $u, r$ coordinates, is the region where $\mu < 1$ and the future light cones have positive expansion: $\partial r/\partial v > 0$, while the future of $\mathcal{A}$ in $\mathcal{Q}$ is the *trapped region* $\mathcal{T}$, where $\mu > 1$ and the future light cones have negative expansion: $\partial r/\partial v < 0$. The future boundary of $\mathcal{T}$ is the *noncentral component* $\mathcal{B} \setminus \mathcal{B}_0$ of the singular boundary $\mathcal{B}$. The function $r$ extends continuously to $\mathcal{B} \setminus \mathcal{B}_0$ where it vanishes. According to the results of [1], $\mathcal{B} \setminus \mathcal{B}_0$ is a strictly spacelike $C^1$ curve, given by

$$(2.5) \qquad \mathcal{B} \setminus \mathcal{B}_0 = \{(u, v^*(u)) \, : \, u \in (\underline{u}^*, 0)\}.$$

Here $v^*$ is a strictly decreasing $C^1$ function in $(\underline{u}^*, 0)$, $v^* > v_0$. Letting

$$(2.6) \qquad \underline{v}^* = \lim_{u \to 0} v^*(u),$$

the central component $\mathcal{B}_0$ of the singular boundary $\mathcal{B}$ is given by

$$(2.7a) \qquad \mathcal{B}_0 = \{(0, v) \, : \, v \in [0, \underline{v}^*]\}$$

if $\mathcal{A}$ is nonempty, and

(2.7b) $$\mathcal{B}_0 = \{(0,v) \,:\, v \in [0,\infty)\}$$

if $\mathcal{A}$ is empty. Also, if $\mathcal{A}$ is nonempty we have:

(2.8) $$\underline{v}_0 = \lim_{u \to 0} v_0(u) \in [0, \underline{v}^*].$$

In [4] we constructed examples where $\mathcal{A}$ is empty and the solutions have a regular extension to $\mathcal{B}_0 \setminus O$, with $r \to \infty$ as $v \to \infty$ on $\mathcal{B}_0$. Then $O$ is a *naked singularity*. We also constructed examples where $\mathcal{B}_0 \setminus O$ is nonempty with $r$ extending continuously to $\mathcal{B}_0$ where it vanishes. Then $\mathcal{B}_0$ corresponds to a *singular future null cone* which has collapsed to a line (see [4]).

Integrating equation (1.20b) yields:

(2.9a) $$\theta_0(t) = e^{\gamma(t)}(\theta_0(0) - I(t))$$

where

(2.9b) $$I(t) = -\int_0^t e^{-\gamma(t')} \zeta_0(t')dt'.$$

The aim of the present section is the proof of the following theorem.

THEOREM 2.1.   *Let $\gamma$ be unbounded. Suppose that either $I$ does not tend to a finite limit as $t \to \infty$ or, otherwise,*

$$\theta_0(0) \neq \lim_{t \to \infty} I(t).$$

*Then $\mathcal{A}$ is nonempty,*

$$\underline{v}_0 = \underline{v}^* = 0,$$

*so that $\mathcal{B}_0 = O$ and both $\mathcal{A}$ and $\mathcal{B} \setminus \mathcal{B}_0$ issue from $O$.*

We recall the following theorem, which was proved in [2].

THEOREM*.   *Let $C_0^+$ be a future light cone with vertex on $\Gamma$ and consider the annular region in $C_0^+$ bounded by two spheres $S_{1,0}$ and $S_{2,0}$, with $S_{2,0}$ in the exterior of $S_{1,0}$. Let $\delta_0$ and $\eta_0$ be the dimensionless size and dimensionless mass content of the region, defined by*:

$$\delta_0 = \frac{r_{2,0}}{r_{1,0}} - 1, \quad \eta_0 = \frac{2(m_{2,0} - m_{1,0})}{r_{2,0}}.$$

*Let $C_1^-$ and $C_2^-$ be the incoming null hypersurfaces through $S_{1,0}$ and $S_{2,0}$ and consider the spheres $S_1$ and $S_2$ at which $C_1^-$ and $C_2^-$ intersect future light cones $C^+$ with vertices on $\Gamma$ in the future of $C_0^+$. Then there are positive constants $c_0 \leq 1/e$ and $c_1 \geq 1$ such that if $\delta_0 \leq c_0$ and*

$$\eta_0 > c_1 \delta_0 \log\left(\frac{1}{\delta_0}\right),$$

then $S_2$ intersects an apparent horizon before $S_1$ reduces to a point on $\Gamma$, that is, there is a future light cone $C_*^+$ such that $S_{2*}$ is a maximal sphere in $C_*^+$ while $r_{1*} > 0$.

In the present context we may apply Theorem* with $C_0^-$ in the role of $C_1^-$, any future light cone $C^+$ intersecting $C_0^-$ in the future of $C_0^+$ in the role of $C_0^+$, and any incoming null hypersurface $C^-$, respecting the spherical symmetry in the exterior of $C_0^-$ in the role of $C_2^-$. Denoting

$$(2.10) \qquad \eta = \frac{2(m - m_0)}{r} = \mu - \mu_0 e^{-s},$$

we conclude that there are positive constants $c_0$ and $c_1$ such that if at some $(t_0, s_0)$ with $t_0 \geq 0$ and $s_0 \in (0, c_0]$ we have

$$\eta(t_0, s_0) > c_1 s_0 \log\left(\frac{1}{s_0}\right),$$

then there is a $t_* \in (t_0, \infty)$ such that the incoming null curve through $(t_0, s_0)$ intersects an apparent horizon at $t = t_*$, so $\kappa \to \infty$ along this curve as $t \to t_*$.

Now if Theorem 2.1 is not true then there is an $\varepsilon > 0$ such that for each $s_0 \in [0, \varepsilon]$ the incoming null curve $s = \chi(t; s_0)$ through $s = s_0$ at $t = 0$ does not intersect an apparent horizon at finite $t$. Let us denote by $\mathcal{R}_\varepsilon$ the region in the half-plane $t \geq 0$,

$$(2.11a) \qquad \mathcal{R}_\varepsilon = \{(t, s) \,:\, t \in [0, \infty) \,\&\, s \in [0, \chi(t; \varepsilon)]\}$$

bounded by the incoming null curves $s = 0$ and $s = \chi(t; \varepsilon)$. Given any positive constant $c$ let us set

$$(2.11b) \qquad \mathcal{R}_\varepsilon^c = \{(t, s) \in \mathcal{R}_\varepsilon \,:\, s \leq c\}.$$

Then according to the above,

$$(2.12) \qquad \eta \leq c_1 s \log\left(\frac{1}{s}\right) \;:\; \text{in } \mathcal{R}_\varepsilon^{c_0}.$$

Let us now consider equation (1.18d) and define:

$$(2.13a) \qquad \psi = e^{-\gamma}(\theta e^s - \theta_0).$$

Using (1.24a) and (1.14) we derive from (1.18d) and (1.20b) the following equation of evolution of $\psi$ along incoming null curves:

$$(2.13b) \qquad \frac{\partial \psi}{\partial t} + \beta \frac{\partial \psi}{\partial s} = \omega \psi + \rho,$$

where

$$(2.13c) \qquad \omega = (1 - \beta)(\kappa - 2) - (\kappa_0 - 2)$$

and

$$(2.13d) \qquad \rho = e^{-\gamma}(\omega \theta_0 + \xi).$$

We shall presently derive estimates for $\omega$ and $\xi$, making use of Lemma 1 as well as the bound (2.12). We begin with an upper bound for $\kappa$. We have:

$$
\begin{aligned}
\kappa(t,s) &= \frac{1}{1-\mu(t,s)} = \frac{1}{1-\mu_0(t)e^{-s}-\eta(t,s)} \\
&\leq \frac{1}{\frac{1}{\kappa_0(t)}-\eta(t,s)} \leq \frac{1}{\frac{e^{-\gamma(t)}}{2\kappa_0(0)}-c_1 s \log\left(\frac{1}{s}\right)}
\end{aligned}
$$

in $\mathcal{R}_\varepsilon^{c_0}$. Thus $(t,s) \in \mathcal{R}_\varepsilon^{c_0}$ and

$$
(2.14) \qquad s \log\left(\frac{1}{s}\right) \leq \frac{e^{-\gamma(t)}}{4c_1\kappa_0(0)}
$$

implies

$$
(2.15) \qquad \kappa(t,s) \leq \kappa_0(0)e^{\gamma(t)}.
$$

Next, we obtain bounds for $\beta$. From (1.21b) and the boundary condition $\nu_0 - \lambda_0 = 0$ (see (1.22)),

$$
(2.16) \qquad (\nu-\lambda)(t,s) = \int_0^s (\kappa(t,s')-1)ds'.
$$

Substituting the estimate (2.15) we obtain, for $(t,s) \in \mathcal{R}_\varepsilon^{c_0}$,

$$
(2.17a) \qquad (\nu-\lambda)(t,s) \leq 4\kappa_0(0)e^{\gamma(t)}s
$$

provided that (2.14) holds. Hence,

$$
(2.17b) \qquad e^{(\nu-\lambda)(t,s)-s} \leq e^{4\kappa_0(0)e^{\gamma(t)}s},
$$

and since for $x \in [0,c]$, $c > 0$, we have

$$
e^x \leq e^c, \quad e^x - 1 \leq \frac{(e^c-1)}{c}x,
$$

while (2.14) implies

$$
(2.18) \qquad 4\kappa_0(0)e^{\gamma(t)}s \leq \frac{1}{c_1},
$$

(recall that $s \leq c_0 \leq 1/e$), also the following hold:

$$
(2.19a) \qquad e^{(\nu-\lambda)(t,s)-s} \leq e^{1/c_1}
$$

$$
(2.19b) \qquad e^{(\nu-\lambda)(t,s)-s} - 1 \leq 4c_1(e^{1/c_1}-1)\kappa_0(0)e^{\gamma(t)}s.
$$

On the other hand for $s \geq 0$, we have

$$
(2.20a) \qquad (\nu-\lambda)(t,s) \geq 0;
$$

hence,

$$
(2.20b) \qquad 1 - e^{(\nu-\lambda)(t,s)-s} \leq 1 - e^{-s} \leq s.
$$

From (2.19a), (2.19b), (2.20b) we conclude, recalling (1.14), that if $(t,s) \in \mathcal{R}_\varepsilon^{c_0}$ and (2.14) holds, then

$$(2.21a) \qquad 0 < 1 - \beta(t,s) \le e^{1/c_1}$$

$$(2.21b) \qquad |\beta(t,s)| \le 4c_1(e^{1/c_1} - 1)\kappa_0(0)e^{\gamma(t)}s.$$

To obtain an estimate for $\omega$ (see (2.13c)), we write:

$$(2.22) \qquad \omega = (1 - \beta)(\kappa - \kappa_0) + \beta(2 - \kappa_0).$$

We have,

$$\kappa - \kappa_0 = \kappa\kappa_0(\mu - \mu_0) = \kappa\kappa_0(\eta - \mu_0(1 - e^{-s})).$$

Using Lemma 1 and estimates (2.12), (2.15), yields:

$$(2.23) \qquad |(\kappa - \kappa_0)(t,s)| \le 8\kappa_0(0)c_1 e^{2\gamma(t)}s\log\left(\frac{1}{s}\right)$$

for $(t,s) \in \mathcal{R}_\varepsilon^{c_0}$ where (2.16) holds. Estimates (2.21a), (2.21b), (2.23), allow us to conclude, in view of the expression (2.22) that:

$$(2.24) \qquad \begin{aligned} |\omega(t,s)| &\le c_2(\kappa_0(0))^2 e^{2\gamma(t)}s\log\left(\frac{1}{s}\right) \\ c_2 &= 16c_1 e^{1/c_1} \end{aligned}$$

for $(t,s) \in \mathcal{R}_\varepsilon^{c_0}$ where (2.14) holds.

To obtain an estimate for $\xi$ we consider the expression (1.24b). By the Schwarz inequality,

$$(2.25) \qquad \xi^2(t,s) \le \int_0^s e^{s'-2\lambda(t,s')}\theta^2(t,s')ds' \cdot \int_0^s e^{-s'+2\nu(t,s')}ds'.$$

Now according to (1.18b),

$$(2.26a) \qquad \frac{\partial(e^s\mu)}{\partial s} = e^{s-2\lambda}\theta^2;$$

thus, in view of the definition (2.10),

$$(2.26b) \qquad e^s\eta(t,s) = \int_0^s e^{s'-2\lambda(t,s')}\theta^2(t,s')ds'.$$

This is the first integral on the right in (2.25). Let us define

$$(2.27a) \qquad \delta(t,s) = e^{s-2\nu(t,s)}\int_0^s e^{-s'+2\nu(t,s')}ds'.$$

By (1.21a), (1.21b),

$$(2.27b) \qquad \frac{\partial\nu}{\partial s} = \frac{1}{2}(e^{2\lambda} - 1 + \theta^2) \ge 0;$$

hence for $s \geq 0$

$$(2.27\text{c}) \qquad \delta(t,s) \leq e^s \int_0^s e^{-s'} ds' = e^s - 1$$

and the second integral on the right in $(2.25)$,

$$(2.27\text{d}) \qquad \int_0^s e^{-s' + 2\nu(t,s')} ds \leq e^{2\nu(t,s)}(1 - e^{-s}).$$

From $(2.25)$, $(2.26\text{b})$, $(2.27\text{d})$, we conclude that

$$(2.28) \qquad \xi^2 \leq \eta e^{2\nu}(e^s - 1).$$

Writing

$$e^{2\nu} = \kappa e^{2(\nu - \lambda)},$$

and using $(2.15)$ and the fact that by $(2.17\text{a})$, $(2.18)$,

$$(2.29\text{a}) \qquad e^{\nu - \lambda} \leq e^{1/c_1},$$

we obtain

$$(2.29\text{b}) \qquad e^{2\nu} \leq 4e^{2/c_1} \kappa_0(0) e^{\gamma(t)}$$

in $\mathcal{R}_\varepsilon^{c_0}$ where $(2.14)$ holds. Hence, using in addition estimate $(2.12)$ we conclude that

$$(2.30\text{a}) \qquad \xi^2 \leq 4c_1 e^{2/c_1} \kappa_0(0) e^{\gamma(t)}(e^s - 1)s \log\left(\frac{1}{s}\right).$$

This also implies that

$$(2.30\text{b}) \qquad |\xi(t,s)| \leq c_2(\kappa_0(0))^2 e^{2\gamma(t)} s \log\left(\frac{1}{s}\right)$$

(see $(2.24)$) for $(t,s) \in \mathcal{R}_\varepsilon^{c_0}$ where $(2.14)$ holds.

We now begin the proof of Theorem 2.1. Since nothing *a priori* is known about the asymptotic behaviour of the integral $I$ defined by $(2.9\text{b})$, we must consider all possibilities. Let

$$(2.31) \qquad l_+ = \limsup_{t \to \infty} I(t) \qquad l_- = \liminf_{t \to \infty} I(t).$$

Then any $l \in [l_-, l_+]$ is a limit value of $I$. The following seven cases exhaust all possibilities:

*Case* 1: $-\infty < l_- = l_+ < \infty$. In this case

$$I(t) \to l \quad \text{as} \quad t \to \infty \quad (l = l_- = l_+)$$

and the hypothesis of Theorem 2.1 states that

$$\theta_0(0) \neq l.$$

*Case* 2: $l_- = l_+ = \infty$.

*Case* 3: $l_- = l_+ = -\infty$.
*Case* 4: $-\infty < l_- < l_+ < \infty$.
*Case* 5: $-\infty < l_- < l_+ = \infty$.
*Case* 6: $-\infty = l_- < l_+ < \infty$.
*Case* 7: $-\infty = l_- < l_+ = \infty$.

In Cases 1 and 4 the integral $I$ is bounded, hence so is $\theta_0 e^{-\gamma}$. Let us set, in these cases,

$$(2.32) \qquad b = \sup_{t \in [0,\infty)} \left| \theta_0(t) e^{-\gamma(t)} \right|.$$

In Case 1, setting

$$(2.33a) \qquad h = |\theta_0(0) - l|$$

we can find a $T > 0$ such that

$$|I_0(t) - l| \leq \frac{2h}{3} \qquad : \ \text{for all } t \geq T;$$

hence,

$$(2.33b) \qquad |\theta_0(t)| \geq \frac{h}{3} e^{\gamma(t)} \qquad : \ \text{for all } t \geq T.$$

In Case 4, setting

$$(2.34a) \qquad h = l_+ - l_-$$

we have

$$\max\{|\theta_0(0) - l_-|, |\theta_0(0) - l_+|\} \geq \frac{h}{2}.$$

It follows that there is an increasing sequence $(t_n : n = 1, 2, ...)$, $t_n \to \infty$ as $n \to \infty$, such that

$$|\theta_0(0) - I(t_n)| \geq \frac{h}{3};$$

hence,

$$(2.34b) \qquad |\theta_0(t_n)| \geq \frac{h}{3} e^{\gamma(t_n)} \qquad : n = 1, 2, \dots.$$

We shall treat Cases 1 and 4 first. From (2.13d) and estimates (2.24) and (2.30b) we have, in view of (2.32),

$$(2.35) \qquad |\rho(t,s)| \leq c_2 (\kappa_0(0))^2 (b+1) e^{2\gamma(t)} s \log\left(\frac{1}{s}\right)$$

in $\mathcal{R}_\varepsilon^{c_0}$ where (2.14) holds. Denoting by

$$\frac{d}{dt} = \frac{\partial}{\partial t} + \beta \frac{\partial}{\partial s},$$

the derivative along incoming null curves with respect to the parameter $t$, we then obtain from (2.13b),

$$(2.36) \qquad \frac{d|\psi|}{dt} \leq c_2(\kappa_0(0))^2 e^{2\gamma(t)} s \log\left(\frac{1}{s}\right)(|\psi| + b + 1)$$

along the incoming null curve

$$(2.37) \qquad s = \chi(t; s_0) \qquad (\chi(0; s_0) = s_0)$$

through $s = s_0 \in [0, \varepsilon]$ at $t = 0$, provided that $s \leq c_0$ and (2.14) holds. Let us define:

$$(2.38) \qquad \varphi(t; s_0) = c_2(\kappa_0(0))^2 \int_0^t e^{2\gamma(t')}\left[s \log\left(\frac{1}{s}\right)\right]_{s=\chi(t'; s_0)} dt'.$$

Integrating (2.36) along the incoming null curve (2.37) then yields:

$$(2.39) \qquad |\psi(t, \chi(t; s_0))| \leq e^{\varphi(t; s_0)}|\psi(0, s_0)| + (b+1)(e^{\varphi(t; s_0)} - 1)$$

provided that

$$(2.40a) \qquad\qquad\qquad \chi(t'; s_0) \leq c_0$$

and (condition 2.14)

$$(2.40b) \qquad \left[s \log\left(\frac{1}{s}\right)\right]_{s=\chi(t'; s_0)} \leq \frac{e^{-\gamma(t')}}{4c_1\kappa_0(0)}$$

hold for all $t' \in [0, t]$.

To proceed we must estimate $\varphi(t; s_0)$. This requires the bound

$$(2.41a) \qquad \kappa(t, s) \leq \frac{3}{2}\kappa_0(t) \;:\; \text{for all} \;\; s \in [0, \chi(t; s_0)].$$

Now, by (2.12),

$$\begin{aligned} \mu(t, s) &= \mu_0(t)e^{-s} + \eta(t, s) \\ &\leq \mu_0(t) + c_1 s \log\left(\frac{1}{s}\right) \end{aligned}$$

for $(t, s) \in \mathcal{R}_\varepsilon^{c_0}$, while (2.41a) is equivalent to:

$$\mu(t, s) \leq \frac{1}{3} + \frac{2}{3}\mu_0(t) \;:\; \text{for all} \;\; s \in [0, \chi(t; s_0)].$$

Consequently, (2.41a) holds when

$$c_1 s \log\left(\frac{1}{s}\right) \leq \frac{1}{3}(1 - \mu_0(t)) = \frac{1}{3\kappa_0(t)}.$$

Thus, by Lemma 1, the bound (2.41a) follows if

$$(2.41b) \qquad\qquad \left[s \log\left(\frac{1}{s}\right)\right]_{s=\chi(t'; s_0)} \leq \frac{e^{-\gamma(t')}}{6c_1\kappa_0(0)}.$$

We note that condition (2.41b) is stronger than condition (2.40b), so only condition (2.41b) together with condition (2.40a) need be considered from this point.

Substituting the bound (2.41a) in (2.16) yields:

$$(2.42) \qquad (\nu - \lambda)(t, \chi(t; s_0)) \leq \left(\frac{3}{2}\kappa_0(t) - 1\right)\chi(t; s_0);$$

hence,

$$(2.43) \qquad \frac{d\chi(t; s_0)}{dt} = \beta(t, \chi(t; s_0)) \geq 1 - e^{\left(\frac{3}{2}\kappa_0(t) - 2\right)\chi(t; s_0)}.$$

Now,

$$f(x) = \begin{cases} (e^x - 1 - x)/x^2 & : \ x \neq 0 \\ 1/2 & : \ x = 0 \end{cases}$$

is a continuous strictly increasing function on the real line. Thus $x \leq 1$ implies

$$f(x) \leq f(1) = e - 2;$$

that is,

$$e^x - 1 \leq x + (e - 2)x^2.$$

Since condition (2.41b) implies that

$$\left(\frac{3}{2}\kappa_0(t) - 2\right)\chi(t; s_0) \leq 1,$$

it follows that

$$1 - e^{\left(\frac{3}{2}\kappa_0 - 2\right)\chi} \geq -\left(\frac{3}{2}\kappa_0 - 2\right)\chi - (e - 2)\left(\frac{3}{2}\kappa_0 - 2\right)^2 \chi^2,$$

and from (2.43),

$$(2.44a) \qquad \frac{d\chi}{dt} \geq -\left(\frac{3}{2}\kappa_0 - 2\right)\chi - (e - 2)\left(\frac{3}{2}\kappa_0 - 2\right)^2 \chi^2,$$

or

$$(2.44b) \qquad \frac{d}{dt}\left(e^{-\frac{3}{2}\gamma + \frac{1}{2}t}\frac{1}{\chi}\right) \leq (e - 2)e^{-\frac{3}{2}\gamma + \frac{1}{2}t}\left(\frac{3}{2}\kappa_0 - 2\right)^2.$$

Let $0 \leq t_0 \leq t_1$. Integrating (2.44b) on $[t_0, t_1]$ yields:

$$(2.44c) \qquad e^{-\frac{3}{2}\gamma(t_1) + \frac{1}{2}t_1}\frac{1}{\chi(t_1)} - e^{-\frac{3}{2}\gamma(t_0) + \frac{1}{2}t_0}\frac{1}{\chi(t_0)} \leq g(t_0, t_1),$$

where

$$(2.44d) \qquad g(t_0, t_1) = (e - 2)\int_{t_0}^{t_1} e^{-\frac{3}{2}\gamma + \frac{1}{2}t}\left(\frac{3}{2}\kappa_0 - 2\right)^2 dt.$$

We have

(2.45a) $$\left(\frac{3}{2}\kappa_0 - 2\right)^2 \leq \frac{9}{4}(\kappa_0 - 1)^2 + \frac{1}{4}$$

and

(2.45b) $$\int_{t_0}^{t_1} e^{-\frac{3}{2}\gamma + \frac{1}{2}t}dt \leq 2e^{-\frac{3}{2}\gamma(t_0)}\left(e^{\frac{1}{2}t_1} - e^{\frac{1}{2}t_0}\right)$$

while, by Lemma 1,

$$\int_{t_0}^{t_1}(\kappa_0 - 1)^2 e^{-\frac{3}{2}\gamma + \frac{1}{2}t}dt = \int_{t_0}^{t_1}(\kappa_0 - 1)\frac{d\gamma}{dt}e^{-\frac{3}{2}\gamma + \frac{1}{2}t}dt$$

$$\leq 2\kappa_0(0)\int_{t_0}^{t_1}\frac{d\gamma}{dt}e^{-\frac{1}{2}\gamma + \frac{1}{2}t}dt$$

$$= 4\kappa_0(0)\left\{e^{-\frac{1}{2}\gamma(t_0) + \frac{1}{2}t_0} - e^{-\frac{1}{2}\gamma(t_1) + \frac{1}{2}t_1} + \frac{1}{2}\int_{t_0}^{t_1}e^{-\frac{1}{2}\gamma + \frac{1}{2}t}dt\right\},$$

which, since

$$\frac{1}{2}\int_{t_0}^{t_1}e^{-\frac{1}{2}\gamma + \frac{1}{2}t}dt \leq e^{-\frac{1}{2}\gamma(t_0)}\left(e^{\frac{1}{2}t_1} - e^{\frac{1}{2}t_0}\right)$$

implies

(2.45c) $$\int_{t_0}^{t_1}(\kappa_0 - 1)^2 e^{-\frac{3}{2}\gamma + \frac{1}{2}t}dt \leq 4\kappa_0(0)e^{\frac{1}{2}t_1}\left(e^{-\frac{1}{2}\gamma(t_0)} - e^{-\frac{1}{2}\gamma(t_1)}\right).$$

By (2.45b) and (2.45c) we conclude, in view of (2.45a), that

(2.46a) $$g(t_0, t_1) \leq (e-2)\left\{\frac{1}{2}e^{-\frac{3}{2}\gamma(t_0)}\left(e^{\frac{1}{2}t_1} - e^{\frac{1}{2}t_0}\right)\right.$$
$$\left. + 9\kappa_0(0)e^{\frac{1}{2}t_1}\left(e^{-\frac{1}{2}\gamma(t_0)} - e^{-\frac{1}{2}\gamma(t_1)}\right)\right\},$$

which implies

(2.46b) $$g(t_0, t_1) \leq c_3\kappa_0(0)e^{\frac{1}{2}t_1} \qquad c_3 = \frac{19}{2}(e-2).$$

Consequently,

(2.47a) $$g(t_0, t_1) \leq \frac{1}{2}e^{-\frac{3}{2}\gamma(t_1) + \frac{1}{2}t_1}\frac{1}{\chi(t_1; s_0)}$$

provided that

(2.47b) $$\chi(t_1; s_0) \leq \frac{e^{-\frac{3}{2}\gamma(t_1)}}{2c_3\kappa_0(0)}.$$

We conclude from (2.44c) that under this condition

$$e^{-\frac{3}{2}\gamma(t_0) + \frac{1}{2}t_0}\frac{1}{\chi(t_0; s_0)} \geq \frac{1}{2}e^{-\frac{3}{2}\gamma(t_1) + \frac{1}{2}t_1}\frac{1}{\chi(t_1; s_0)};$$

that is,

$$(2.48) \qquad \chi(t_0; s_0) \le 2e^{\frac{3}{2}(\gamma(t_1) - \gamma(t_0)) - \frac{1}{2}(t_1 - t_0)} \chi(t_1; s_0)$$

holds.

We turn to the estimation of $\varphi(t; s_0)$ defined by (2.38). Setting $t_0 = t'$, $t_1 = t$ in (2.48) and noting that for $x > 0$ we have

$$x^{1/2} \log\left(\frac{1}{x}\right) \le \frac{2}{e},$$

and we obtain

$$(2.49a) \quad \chi(t'; s_0) \log\left(\frac{1}{\chi(t'; s_0)}\right) \quad \le \quad \frac{2}{e}(\chi(t'; s_0))^{1/2}$$

$$\le \quad \frac{2^{3/2}}{e} e^{\frac{3}{4}(\gamma(t) - \gamma(t')) - \frac{1}{4}(t - t')}(\chi(t; s_0))^{1/2}.$$

Hence,

$$(2.49b)$$

$$\int_0^t e^{2\gamma(t')} \chi(t'; s_0) \log\left(\frac{1}{\chi(t'; s_0)}\right) dt' \quad \le \quad \frac{2^{3/2}}{e} e^{2\gamma(t)} \int_0^t e^{-\frac{1}{4}(t - t')} dt' \cdot (\chi(t; s_0))^{1/2}$$

$$\le \quad \frac{2^{7/2}}{e} e^{2\gamma(t)}(\chi(t; s_0))^{1/2},$$

and

$$(2.50) \qquad \varphi(t; s_0) \quad \le \quad c_4 (\kappa_0(0))^2 e^{2\gamma(t)}(\chi(t; s_0))^{1/2}$$

$$c_4 = \frac{2^{7/2}}{e} c_2.$$

This estimate holds provided that condition (2.40a), condition (2.41b) with $t$ replaced by $t'$, that is,

$$(2.51a) \qquad \chi(t'; s_0) \log\left(\frac{1}{\chi(t'; s_0)}\right) \le \frac{e^{-\gamma(t)}}{6c_1 \kappa_0(0)},$$

is satisfied for all $t' \in [0, t]$ and condition (2.47b) with $t_1$ replaced by $t$, that is,

$$(2.51b) \qquad \chi(t; s_0) \le \frac{e^{-\frac{3}{2}\gamma(t)}}{2c_3 \kappa_0(0)},$$

is satisfied as well.

By virtue of estimate (2.50), if

$$(2.52) \qquad \chi(t; s_0) \le \left(\min\left\{\log 2, \frac{h}{48(b+1)}\right\}\right)^2 \frac{e^{-4\gamma(t)}}{c_4^2 (\kappa_0(0))^4}$$

then we have

$$(2.53a) \qquad e^{\varphi(t; s_0)} \le 2,$$

and

(2.53b) $$e^{\varphi(t;s_0)} - 1 \leq 2\varphi(t;s_0) \leq \frac{h}{24(b+1)}.$$

Therefore, if also

(2.54a) $$\sup_{s \in [0,s_0]} |\psi(0;s)| \leq \frac{h}{48}$$

then from (2.39) we conclude that

(2.54b) $$|\psi(t, \chi(t;s_0))| \leq \frac{h}{12}.$$

By (2.13a) and (2.34b) this implies:

(2.55a) $$|\theta(t_n, s)e^s| \geq \frac{h}{4}e^{\gamma(t_n)} \qquad : \ n = 1, 2, \ldots$$

for all $s \in [0, s_n]$, where

(2.55b) $$s_n = \chi(t_n; s_0).$$

For, if conditions (2.40a), (2.51a), (2.51b), (2.52), (2.54a) hold for $s_0$, then they hold *a fortiori* if $s_0$ is replaced by $s_0' \in [0, s_n]$; thus (2.54b) also holds with $s_0$ replaced by $s_0'$. Since for each $s \in [0, s_n]$ there is a $s_0' \in [0, s_0]$ such that $\chi(t_n; s_0') = s$, the result follows. We remark that in Case 1, the sequence $(t_n : n = 1, 2, \ldots)$ can be chosen to be an arbitrary increasing sequence contained in $[T, \infty)$ such that $t_n \to \infty$ as $n \to \infty$. Now, (2.55a) implies (see (2.26b)):

(2.56a) $$\begin{aligned} \eta(t_n, s_n) &= e^{-s_n} \int_0^{s_n} e^s \frac{\theta^2(t_n, s)}{\kappa(t_n, s)} ds \\ &\geq \frac{h^2 e^{\gamma(t_n)}}{64 \kappa_0(0)} (1 - e^{-s_n}) \end{aligned}$$

where we have used the bound (2.15). Since $s_n \leq c_0$ and the function

$$\frac{1 - e^{-x}}{x}$$

is decreasing for $x > 0$, we have

$$\frac{1 - e^{-s_n}}{s_n} \geq \frac{1 - e^{-c_0}}{c_0}.$$

Thus (2.56a) implies:

(2.56b) $$\begin{aligned} \eta(t_n, s_n) &\geq \frac{c_5 h^2}{\kappa_0(0)} e^{\gamma(t_n)} s_n \\ c_5 &= \frac{1}{64} \frac{1 - e^{-c_0}}{c_0}. \end{aligned}$$

The lower bound (2.56b) contadicts the upper bound (2.12) if we choose

$$(2.57) \qquad s_n = e^{-c_6 h^2 e^{\gamma(t_n)}/\kappa_0(0)}$$

where $c_6$ is some constant such that

$$c_6 > \frac{c_5}{c_1}.$$

We shall show that the choice (2.57) is, for sufficiently large $n$, consistent with conditions (2.40a), (2.51a), (2.51b), (2.52), (2.54a).

We note that in the above conditions $t$ stands for $t_n$. Thus conditions (2.51b) and (2.52) read:

$$(2.58a) \qquad s_n \leq \frac{e^{-\frac{3}{2}\gamma(t_n)}}{2c_3\kappa_0(0)}$$

and

$$(2.58b) \qquad s_n \leq \left(\min\left\{\log 2, \frac{h}{48(b+1)}\right\}\right)^2 \frac{e^{-4\gamma(t_n)}}{c_4^2(\kappa_0(0))^4},$$

respectively. Since, according to the hypotheses of Theorem 2.1,

$$\gamma(t_n) \to \infty \ \text{as} \ n \to \infty,$$

both these conditions are satisfied by the choice (2.57) if $n$ is sufficiently large. Also, conditions (2.40a) and (2.51a) are satisfied at $t' = t_n$ if $n$ is sufficiently large:

$$(2.59a) \qquad s_n \ \leq \ c_0$$
$$(2.59b) \qquad s_n \log\left(\frac{1}{s_n}\right) \ \leq \ \frac{e^{-\gamma(t_n)}}{6c_1\kappa_0(0)}.$$

Let $t_*$ be the least value of $t \in [0, t_n]$ such that (2.40a) and (2.51a) are satisfied for all $t' \in [t, t_n]$. Then if $t_* > 0$, either

$$(2.60a) \qquad \chi(t_*; s_0) = c_0$$

or

$$(2.60b) \qquad \chi(t_*; s_0) \log\left(\frac{1}{\chi(t_*; s_0)}\right) = \frac{e^{-\gamma(t_*)}}{6c_1\kappa_0(0)}.$$

However under these circumstances estimate (2.48) holds with $t_0$, $t_1$ replaced by $t_*$, $t_n$, respectively

$$\chi(t_*; s_0) \leq 2e^{\frac{3}{2}(\gamma(t_n) - \gamma(t_*)) - \frac{1}{2}(t_n - t_*)} s_n;$$

so, *a fortiori*,

$$(2.61) \qquad \chi(t_*; s_0) \leq 2e^{\frac{3}{2}\gamma(t_n)} s_n.$$

We see that the choice (2.57) contradicts both (2.60a), (2.60b), if $n$ is sufficiently large. We conclude that $t_* = 0$ so conditions (2.40a) and (2.51a) are satisfied for all $t' \in [0, t_n]$ and (2.61) reads

$$(2.62) \qquad s_0 \leq 2e^{\frac{3}{2}\gamma(t_n)}s_n \to 0 \ \text{ as } n \to \infty.$$

In view of the fact that solutions of bounded variation have the property that $\theta(t, s)$ is at each $t$, in particular at $t = 0$, continuous from the right with respect to $s$, condition (2.54a) also follows for sufficiently large $n$. We conclude that the choice (2.57) is, for sufficiently large $n$, consistent with all conditions. We therefore reach a contradiction if we suppose Theorem 2.1 to be false in Cases 1 and 4.

To treat the remaining Cases 2, 3, 5, 6, 7, in which $I$ is unbounded, we set:

$$(2.63a) \qquad b(t) = \sup_{t' \in [0,t]} \left| \theta_0(t')e^{-\gamma(t')} \right|.$$

Then $b$ is a nondecreasing function tending to infinity as $t \to \infty$, and there is an increasing sequence $(t_n : n = 1, 2, ...)$, $t_n \to \infty$ as $n \to \infty$, such that

$$(2.63b) \qquad \left| \theta_0(t_n)e^{-\gamma(t_n)} \right| = b_0(t_n) \to \infty \ \text{ as } n \to \infty.$$

In $[0, t_n]$ we have, as in (2.36),

$$\frac{d|\psi|}{dt} \leq c_2(\kappa_0(0))^2 e^{2\gamma(t)} s \log\left(\frac{1}{s}\right)(|\psi| + b(t_n) + 1);$$

hence integrating we obtain

$$(2.64) \qquad |\psi(t_n, s_n)| \leq e^{\varphi(t_n; s_0)}|\psi(0, s_0)| + (b(t_n) + 1)\left(e^{\varphi(t_n; s_0)} - 1\right)$$

where again

$$s_n = \chi(t_n; s_0).$$

As before, the conditions

$$(2.65a) \qquad \chi(t; s_0) \quad \leq \quad c_0 \ : \ \text{for all } t \in [0, t_n],$$

$$(2.65b) \quad \chi(t; s_0)\log\left(\frac{1}{\chi(t; s_0)}\right) \quad \leq \quad \frac{e^{-\gamma(t)}}{6c_1\kappa_0(0)} \ : \ \text{for all } t \in [0, t_n],$$

and

$$(2.65c) \qquad s_n \leq \frac{e^{-\frac{3}{2}\gamma(t_n)}}{2c_3\kappa_0(0)}$$

imply

$$(2.66) \qquad \varphi(t_n; s_0) \leq c_4(\kappa_0(0))^2 e^{2\gamma(t_n)}s_n^{1/2}$$

(see (2.50)). Thus if

$$(2.67) \qquad s_n \leq \left( \min\left\{ \log 2, \frac{1}{4(b(t_n)+1)} \right\} \right)^2 \frac{e^{-4\gamma(t_n)}}{c_4^2(\kappa_0(0))^4},$$

then we have

$$(2.68a) \qquad e^{\varphi(t_n; s_0)} \leq 2$$

and

$$(2.68b) \qquad e^{\varphi(t_n; s_0)} - 1 \leq 2\varphi(t_n; s_0) \leq \frac{1}{2(b(t_n)+1)}.$$

Therefore if also

$$(2.69) \qquad \sup_{s \in [0, s_0]} |\psi(0, s)| \leq \frac{1}{4}$$

we conclude from (2.64) that

$$(2.70a) \qquad |\psi(t_n, s_n)| \leq 1.$$

In fact, since conditions (2.65a)–(2.65c), (2.67), (2.69), hold *a fortiori* if $s_0$ is replaced by $s_0' \in [0, s_0]$, we conclude that

$$(2.70b) \qquad \sup_{s \in [0, s_n]} |\psi(t_n, s)| \leq 1.$$

Hence,

$$(2.71) \qquad \begin{aligned} |\theta(t_n, s)e^s| &\geq |\theta_0(t_n)| - e^{\gamma(t_n)}|\psi(t_n, s)| \\ &\geq (b(t_n) - 1)e^{\gamma(t_n)} \geq \frac{1}{2}b(t_n)e^{\gamma(t_n)} \end{aligned}$$

for all $s \in [0, s_n]$, if $n$ is large enough. It follows that (see (2.56a), (2.56b))

$$(2.72) \qquad \begin{aligned} \eta(t_n, s_n) &\geq c_7 \frac{(b(t_n))^2}{\kappa_0(0)} e^{\gamma(t_n)} s_n \\ c_7 &= \frac{1}{16} \frac{1 - e^{-1/c_0}}{1/c_0}. \end{aligned}$$

The lower bound (2.72) contradicts the upper bound (2.12) if we choose

$$(2.73) \qquad s_n = e^{-c_8(b(t_n))^2 e^{\gamma(t_n)}/\kappa_0(0)},$$

where $c_8$ is some constant such that

$$c_8 > \frac{c_7}{c_1}.$$

Since

$$b(t_n) \to \infty \quad \text{as} \quad n \to \infty,$$

as well as

$$\gamma(t_n) \to \infty \quad \text{as} \quad n \to \infty,$$

the conditions (2.65c), (2.67) are satisfied by the choice (2.73) if $n$ is sufficiently large. Also, the conditions (2.65a), (2.65b) are satisfied at $t = t_n$. The continuity argument which we applied previously then shows that conditions (2.65a), (2.65b) are satisfied for all $t \in [0, t_n]$ and, moreover, (2.62) holds, which implies that condition (2.69) is verified as well. We conclude that the choice (2.73) is, for sufficiently large $n$, consistent with all conditions. We therefore again reach a contradiction if we suppose Theorem 2.1 to be false in the remaining Cases 2, 3, 5, 6, 7. This completes the proof of Theorem 2.1.

## 3. The second instability theorem

In the following we confine attention to the case not covered by Theorem 2.1, the case where $I$ tends to a finite limit as $t \to \infty$ and

$$(3.1a) \qquad \theta_0(0) = \lim_{t \to \infty} I(t).$$

Following the proof of Theorem 2.1 we see that the argument of Cases 1 and 4 still applies if there is a positive constant $p < 1$ such that

$$\limsup_{t \to \infty} \left\{ |\theta_0(0) - I(t)| e^{\frac{1}{2}p\gamma(t)} \right\} \neq 0.$$

Therefore we can assume in the following that

$$(3.1b) \qquad (\theta_0(0) - I(t)) e^{\frac{1}{2}p\gamma(t)} \to 0 \quad \text{as} \quad t \to \infty$$

for all positive constants $p < 1$.

Let us define a function $\tau$ by

$$(3.2a) \qquad \frac{\partial \tau}{\partial t} + \beta \frac{\partial \tau}{\partial s} = \omega, \qquad \tau(0, s) = 0.$$

We then have

$$(3.2b) \qquad \tau(t, \chi(t; s_0)) = \int_0^t \omega(t', \chi(t'; s_0)) dt'.$$

Let us also define the functions

$$(3.3a) \qquad \tilde{\psi} = e^{-\tau} \psi$$
$$(3.3b) \qquad \tilde{\rho} = e^{-\tau} \rho.$$

Then, from (2.13b) we have

$$(3.3c) \qquad \frac{\partial \tilde{\psi}}{\partial t} + \beta \frac{\partial \tilde{\psi}}{\partial s} = \tilde{\rho}, \qquad \tilde{\psi}(0, s) = \psi(0, s).$$

Thus, if we also define a function $\sigma$ by

$$(3.3\text{d}) \qquad \frac{\partial \sigma}{\partial t} + \beta \frac{\partial \sigma}{\partial s} = \tilde{\rho}, \qquad \sigma(0,s) = 0,$$

that is,

$$(3.3\text{e}) \qquad \sigma(t, \chi(t; s_0)) = \int_0^t \tilde{\rho}(t', \chi(t'; s_0)) dt',$$

then

$$(3.3\text{f}) \qquad \tilde{\psi} = \hat{\psi} + \sigma$$

where $\hat{\psi}$ satisfies

$$(3.3\text{g}) \qquad \frac{\partial \hat{\psi}}{\partial t} + \beta \frac{\partial \hat{\psi}}{\partial s} = 0, \qquad \hat{\psi}(0,s) = \psi(0,s)$$

so that

$$(3.3\text{h}) \qquad \hat{\psi}(t, \chi(t; s_0)) = \psi(0, s_0).$$

Let us fix

$$p = \frac{1}{2}$$

in (3.1b). Then from (2.9a) there is a constant $b$ such that

$$(3.4) \qquad |\theta_0(t)| \leq b e^{\frac{3}{4}\gamma(t)} \quad : \text{for all } t \geq 0.$$

If the conclusion of Theorem 2.1 is false, there is an $\varepsilon > 0$ such that for each $s_0 \in [0, \varepsilon]$ the incoming null curve $s = \chi(t; s_0)$ through $s = s_0$ at $t = 0$ does not intersect an apparent horizon at finite $t$. Estimates (2.24) and (2.30b) then hold in $\mathcal{R}_\varepsilon^{c_0}$ where condition (2.14) holds. It follows that:

$$(3.5) \qquad |\rho(t, s)| \leq c_2(\kappa_0(0))^2(b+1) e^{\frac{7}{4}\gamma(t)} s \log\left(\frac{1}{s}\right)$$

in $\mathcal{R}_\varepsilon^{c_0}$ where (2.14) holds. Also, conditions (2.51a) and (2.51b) imply (2.48) with $t_1, t_0$ replaced by $t, t' \in [0, t]$, respectively, that is

$$(3.6) \qquad s' \leq 2 e^{\frac{3}{2}(\gamma(t) - \gamma(t')) - \frac{1}{4}(t - t')} s$$

with

$$s' = \chi(t'; s_0), \quad s = \chi(t; s_0).$$

The same conditions imply estimate (2.50), which if

$$(3.7) \qquad \chi(t; s_0) \leq (\log 2)^2 \frac{e^{-4\gamma(t)}}{c_4^2(\kappa_0(0))^4}$$

then it, in turn, implies

$$(3.8) \qquad e^{\varphi(t; s_0)} \leq 2.$$

Hence, in view of the fact that by estimate (2.24),

$$(3.9) \qquad\qquad |\tau(t, \chi(t; s_0))| \leq \varphi(t; s_0)$$

we have

$$(3.10a) \qquad
\begin{aligned}
|\tilde{\rho}(t, s)| &\leq 2|\rho(t, s)| \\
&\leq 2c_2(\kappa_0(0))^2 (b+1) e^{\frac{7}{4}\gamma(t)} s \log\left(\frac{1}{s}\right).
\end{aligned}$$

Moreover, if (2.51b) and (3.7) hold with $t$ replaced by $t'$, for all $t' \in [0, t]$, then also (3.8) and (3.10a) hold with $t, s = \chi(t; s_0)$, replaced by $t', s' = \chi(t'; s_0)$, for all $t' \in [0, t]$. It follows that:

$$(3.10b) \qquad
\begin{aligned}
|\sigma(t, s)| &\leq \int_0^t |\tilde{\rho}(t', s')| dt' \\
&\leq 2c_2(\kappa_0(0))^2 (b+1) \int_0^t e^{\frac{7}{4}\gamma(t')} s' \log\left(\frac{1}{s'}\right) dt'.
\end{aligned}$$

Now, by virtue of (3.6) we have (see (2.49a))

$$
\begin{aligned}
s' \log\left(\frac{1}{s'}\right) &\leq \frac{2}{e} s'^{1/2} \\
&\leq \frac{2^{3/2}}{e} e^{\frac{3}{4}(\gamma(t) - \gamma(t')) - \frac{1}{4}(t - t')} s^{1/2};
\end{aligned}$$

hence

$$
\begin{aligned}
&\int_0^t e^{\frac{7}{4}\gamma(t')} s' \log\left(\frac{1}{s'}\right) dt' \\
&\leq \frac{2^{3/2}}{e} e^{\frac{7}{4}\gamma(t)} \int_0^t e^{-\frac{1}{4}(t - t')} dt' \cdot s^{1/2} \\
&\leq \frac{2^{7/2}}{e} e^{\frac{7}{4}\gamma(t)} s^{1/2}.
\end{aligned}$$

Thus we obtain

$$(3.10c) \qquad\qquad |\sigma(t, s)| \leq 2c_4(\kappa_0(0))^2 e^{\frac{7}{4}\gamma(t)} s^{1/2}$$

(see (2.50)), which, being valid for all $s \in [0, s_*]$ implies

$$(3.10d) \qquad
\begin{aligned}
||\sigma(t)||_{L^2(0, s_*)} &:= \sqrt{\int_0^{s_*} \sigma^2(t, s) ds} \\
&\leq c_4(\kappa_0(0))^2 e^{\frac{7}{4}\gamma(t)} s_*.
\end{aligned}$$

On the other hand, from (3.3h),

$$(3.11) \qquad
\begin{aligned}
||\hat{\psi}(t)||^2_{L^2(0, s_*)} &:= \int_0^{s_*} \hat{\psi}^2(t, s) ds \\
&= \int_0^{s_{0*}} \psi^2(0, s_0) \frac{\partial \chi}{\partial s_0}(t; s_0) ds_0,
\end{aligned}$$

where
$$s = \chi(t; s_0), \quad s_* = \chi(t; s_{0*}).$$

Now, differentiating the equation
$$\frac{d\chi(t; s_0)}{dt} = \beta(t, \chi(t; s_0)), \quad \chi(0; s_0) = s_0$$

with respect to $s_0$ yields:

(3.12a) $$\frac{d}{dt}\left(\frac{\partial\chi}{\partial s_0}(t; s_0)\right) = \frac{\partial\beta}{\partial s}(t, \chi(t; s_0))\frac{\partial\chi}{\partial s_0}(t; s_0), \quad \frac{\partial\chi}{\partial s_0}(0; s_0) = 1.$$

By equations (1.18c) and (2.13c),

(3.12b) $$\frac{\partial\beta}{\partial s} = 2 - \kappa_0 - \omega.$$

Thus, substituting in (3.12a) and integrating,
$$\frac{\partial\chi}{\partial s_0}(t; s_0) = \exp\left[\int_0^t (2 - \kappa_0(t') - \omega(t', \chi(t'; s_0)))dt'\right]$$

or, in view of (2.1) and (3.2b), we see that

(3.12c) $$\frac{\partial\chi}{\partial s_0}(t; s_0) = e^{t - \gamma(t) - \tau(t, \chi(t; s_0))}.$$

By (3.8), (3.9) we then have

(3.12d) $$\frac{\partial\chi}{\partial s_0}(t; s_0) \geq \frac{1}{2}e^{t - \gamma(t)}.$$

Substituting in (3.11) we conclude that

(3.13) $$||\hat{\psi}(t)||_{L^2(0, s_*)} \geq \frac{1}{2^{1/2}}e^{\frac{1}{2}t - \frac{1}{2}\gamma(t)}s_{0*}^{1/2}h(s_{0*}),$$

where

(3.14) $$h(s) = \sqrt{\frac{1}{s}\int_0^s \psi^2(0, s')ds'}.$$

Now we wish to achieve

(3.15a) $$||\hat{\psi}(t)||_{L^2(0, s_*)} \geq 2||\sigma||_{L^2(0, s_*)}$$

so that

(3.15b) $$\begin{aligned}||\tilde{\psi}(t)||_{L^2(0, s_*)} &\geq ||\hat{\psi}(t)||_{L^2(0, s_*)} - ||\sigma||_{L^2(0, s_*)} \\ &\geq \frac{1}{2}||\hat{\psi}(t)||_{L^2(0, s_*)}^2\end{aligned}$$

(see (3.3f)). In view of estimates (3.10d) and (3.13), this requires an upper bound for $s_*$ in terms of $s_{0*}$, which in turn requires an upper bound for $\beta$. Since $\eta \geq 0$, we have
$$\kappa(t, s) = \frac{1}{1 - \mu(t, s)} \geq \frac{1}{1 - \mu_0(t)e^{-s}};$$

hence (see (2.16)),

$$(\nu - \lambda)(t,s) \geq \int_0^s \left( \frac{1}{1 - \mu_0(t)e^{-s'}} - 1 \right) ds'$$
$$= \log \left( \frac{1 - \mu_0(t)e^{-s}}{1 - \mu_0(t)} \right),$$

and

(3.16a)
$$\beta(t,s) = 1 - e^{(\nu-\lambda)(t,s)-s}$$
$$\leq 1 - e^{-2s} - (e^{-s} - e^{-2s})\kappa_0(t).$$

Since

$$1 - e^{-2s} \leq 2s - \frac{3}{2}s^2 \ : \text{if} \ s \leq c_0,$$

(3.16a) implies

(3.16b)
$$\beta(t,s) \leq (2 - \kappa_0(t))s + \frac{3}{2}(\kappa_0(t) - 1)s^2.$$

Thus,

(3.17a)
$$\frac{d\chi}{dt} \leq (2 - \kappa_0)\chi + \frac{3}{2}(\kappa_0 - 1)\chi^2$$

or

(3.17b)
$$\frac{d}{dt}\left( e^{t-\gamma}\frac{1}{\chi} \right) \geq -\frac{3}{2}(\kappa_0 - 1)e^{t-\gamma}.$$

Integrating (3.17b) on $[0,t]$ yields

(3.17c)
$$\frac{e^{t-\gamma(t)}}{s} - \frac{1}{s_0} \geq -\frac{3}{2}\int_0^t (\kappa_0(t') - 1)e^{t'-\gamma(t')}dt'$$
$$= -\frac{3}{2}\int_0^t \frac{d\gamma}{dt'}e^{t'-\gamma(t')}dt'$$
$$= -\frac{3}{2}\left( 1 - e^{t-\gamma(t)} + \int_0^t e^{t'-\gamma(t')}dt' \right)$$
$$\geq -\frac{3}{2}e^t(1 - e^{-\gamma(t)}) \geq -\frac{3}{2}e^t,$$

where $s = \chi(t; s_0)$. It follows that if

(3.18a)
$$s_0 \leq \frac{1}{3}e^{-t},$$

then

(3.18b)
$$s \leq 2e^{t-\gamma(t)}s_0.$$

In particular, if

(3.19a)
$$s_{0*} \leq \frac{1}{3}e^{-t},$$

then

$$(3.19b) \qquad s_* \leq 2e^{t-\gamma(t)}s_{0*}.$$

Substituting this in (3.10d) we conclude that under condition (3.19a),

$$(3.20) \qquad h(s_{0*}) \geq 2^{5/2}c_4(\kappa_0(0))^2 e^{\frac{1}{2}t+\frac{5}{4}\gamma(t)}s_{0*}^{1/2}$$

implies (3.15a), hence also (3.15b).

Now, we have (see (2.56a)):

$$(3.21) \qquad \begin{aligned} \eta(t,s_*) &= e^{-s_*}\int_0^{s_*} e^s \frac{\theta^2(t,s)}{\kappa(t,s)}ds \\ &\geq \frac{e^{-2c_0}}{2\kappa_0(0)}e^{-\gamma(t)}\int_0^{s_*} e^{2s}\theta^2(t,s)ds \end{aligned}$$

where we have used the bound (2.15) and the condition

$$s_* \leq c_0.$$

From (2.13a) and (3.3a), we have

$$(3.22a) \qquad e^s\theta = \theta_0 + e^{\gamma+\tau}\tilde{\psi}.$$

It follows that

$$(3.22b) \qquad ||e^s\theta(t)||_{L^2(0,s_*)} \geq ||e^{\gamma+\tau}\tilde{\psi}(t)||_{L^2(0,s_*)} - ||\theta_0(t)||_{L^2(0,s_*)}.$$

By (3.8) and (3.15b),

$$(3.23a) \qquad ||e^{\gamma+\tau}\tilde{\psi}(t)||_{L^2(0,s_*)} \geq \frac{1}{4}e^{\gamma(t)}||\hat{\psi}(t)||_{L^2(0,s_*)}$$

while by (3.4) and (3.19b),

$$(3.23b) \qquad ||\theta_0(t)||_{L^2(0,s_*)} = |\theta_0(t)|s_*^{1/2} \leq 2^{1/2}e^{\frac{1}{2}t+\frac{1}{4}\gamma(t)}bs_{0*}^{1/2}.$$

In view of the lower bound (3.13) we conclude that

$$(3.24a) \qquad h(s_{0*}) \geq 16be^{-\frac{1}{4}\gamma(t)}$$

implies

$$(3.24b) \qquad ||\theta_0(t)||_{L^2(0,s_*)} \leq \frac{1}{8}e^{\gamma(t)}||\hat{\psi}||_{L^2(0,s_*)},$$

and

$$(3.24c) \qquad \begin{aligned} ||e^s\theta(t)||_{L^2(0,s_*)} &\geq \frac{1}{8}e^{\gamma(t)}||\hat{\psi}(t)||_{L^2(0,s_*)} \\ &\geq \frac{1}{2^{7/2}}e^{\frac{1}{2}t+\frac{1}{2}\gamma(t)}s_{0*}^{1/2}h(s_{0*}). \end{aligned}$$

Substituting this in (3.21) yields:

$$(3.25) \qquad \eta(t, s_*) \geq \frac{e^{-2c_0}}{2^8 \kappa_0(0)} e^t s_{0*} h^2(s_{0*}).$$

On the other hand, the bound (2.12) at $(t, s_*)$ reads

$$(3.26a) \qquad \eta(t, s_*) \leq c_1 s_* \log\left(\frac{1}{s_*}\right),$$

which by (3.19b) implies

$$(3.26b) \qquad \eta(t, s_*) \leq 2c_1 e^{t - \gamma(t)} s_{0*} \left\{ \log\left(\frac{1}{s_{0*}}\right) + \gamma(t) - t \right\}.$$

This contradicts (3.25) if

$$(3.27) \qquad \begin{aligned} h^2(s_{0*}) &> c_9 \kappa_0(0) e^{-\gamma(t)} \left\{ \log\left(\frac{1}{s_{0*}}\right) + \gamma(t) - t \right\} \\ c_9 &= 2^9 e^{2c_0} c_1. \end{aligned}$$

Summarizing, the hypothesis that the conclusion of Theorem 2.1 is false leads to a contradiction if for every $\varepsilon > 0$ there is a $s_{0*} \in (0, \varepsilon]$ and a $t \in [0, \infty)$ such that the requirements (3.20), (3.24a) and (3.27) are satisfied, and, moreover, the conditions (2.40a), (2.51a), (2.51b), (3.7), that is, with

$$s'_* = \chi(t'; s_{0*}),$$

the conditions:

$$(3.28a) \qquad s'_* \leq c_0 \quad : \text{for all } t' \in [0, t]$$

$$(3.28b) \qquad s'_* \log\left(\frac{1}{s'_*}\right) \leq \frac{e^{-\gamma(t')}}{6c_1 \kappa_0(0)} \quad : \text{for all } t' \in [0, t]$$

$$(3.28c) \qquad s'_* \leq \frac{e^{-\frac{3}{2}\gamma(t')}}{2c_3 \kappa_0(0)} \quad : \text{for all } t' \in [0, t]$$

$$(3.28d) \qquad s'_* \leq (\log 2)^2 \frac{e^{-4\gamma(t')}}{c_4^2 (\kappa_0(0))^4} \quad : \text{for all } t' \in [0, t],$$

respectively, hold, and, finally, condition (3.19a) holds as well.

Given $s_{0*}$, let us define $t$ by

$$(3.29a) \qquad t + 5\gamma(t) = \log\left(\frac{1}{s_{0*}}\right),$$

i.e.

$$(3.29b) \qquad s_{0*} = e^{-t - 5\gamma(t)}.$$

Then if $t$ is large enough so that

$$(3.30) \qquad \gamma(t) \geq \frac{1}{5} \log 3,$$

then condition (3.19a) is verified. This implies (3.19b) with $t$ replaced by any $t' \in [0, t]$, that is:

$$(3.31a) \qquad s'_* \leq 2e^{t' - \gamma(t')} s_{0*} \;\; : \text{for all } t' \in [0, t].$$

Substituting (3.29b) we obtain

$$(3.31b) \qquad s'_* \leq 2e^{-5\gamma(t)} \;\; : \text{for all } t' \in [0, t].$$

It follows that (3.28a), (3.28c), (3.28d) are verified if $t$ is large enough so that

$$(3.32a) \qquad \gamma(t) \;\; \geq \;\; \frac{1}{5} \log \left( \frac{2}{c_0} \right)$$

$$(3.32b) \qquad \gamma(t) \;\; \geq \;\; \frac{2}{7} \log(4c_3 \kappa_0(0))$$

$$(3.32c) \qquad \gamma(t) \;\; \geq \;\; 2 \log \left( \frac{c_4(\kappa_0(0))^2}{\log 2} \right) + \log 2,$$

respectively. Also, since

$$s'_* \log \left( \frac{1}{s'_*} \right) \leq \frac{2}{e} s'^{1/2}_*,$$

(3.28b) holds *a fortiori* if

$$(3.33a) \qquad s'_* \leq \frac{e^2 e^{-2\gamma(t')}}{144 c_1^2 (\kappa_0(0))^2} \;\; : \text{for all } t' \in [0, t]$$

which by virtue of (3.31b) is verified if $t$ is large enough so that

$$(3.33b) \qquad \gamma(t) \geq \frac{1}{3} \left\{ 2 \log(12 c_1 \kappa_0(0)) + \log 2 - 2 \right\}.$$

We turn to the requirements (3.20), (3.24a) and (3.27). Substituting the definition (3.29b), these requirements read

$$(3.34a) \qquad h(e^{-t - 5\gamma(t)}) \;\; \geq \;\; 2^{5/2} c_4 (\kappa_0(0))^2 e^{-\frac{5}{4}\gamma(t)}$$

$$(3.34b) \qquad h(e^{-t - 5\gamma(t)}) \;\; \geq \;\; 16 b e^{-\frac{1}{4}\gamma(t)}$$

$$(3.34c) \qquad h(e^{-t - 5\gamma(t)}) \;\; > \;\; (6c_9)^{1/2} (\kappa_0(0))^{1/2} (\gamma(t))^{1/2} e^{-\frac{1}{2}\gamma(t)},$$

respectively. For large $t$, (3.34b) is the strongest requirement. Let us assume that

$$(3.35) \qquad \limsup_{t \to \infty} \left\{ e^{\frac{1}{4}\gamma(t)} h(e^{-t - 5\gamma(t)}) \right\} = \infty.$$

Then given any $T \in (0, \infty)$ we can find a $t \in [T, \infty)$ such that with $s_{0*}$ defined by (3.29b) all conditions and requirements are satisfied. We have therefore proved:

THEOREM 3.1.   *Let $\gamma$ be unbounded and let*

$$I(t) \to \theta_0(0) \;\; \text{as } t \to \infty.$$

*Let $g(s)$ be the function on $(0, 1)$ defined by*

$$g(s) = e^{-\frac{1}{4}\gamma(t)}, \quad s = e^{-t-5\gamma(t)}.$$

*Consider the function*

$$h(s) = \sqrt{\frac{1}{s} \int_0^s \left( e^{s'} \theta(0, s') - \theta(0, 0) \right)^2 ds'}$$

*defined for $s > 0$. Then if*

$$\limsup_{s \to 0+} \left\{ \frac{h(s)}{g(s)} \right\} = \infty,$$

*then the conclusion of Theorem 2.1 holds.*

## 4. The exceptional set

We now investigate the subset $\mathcal{E}$ of the space of initial conditions of bounded variation on $C_0^+$, which lead to the formation of a singular boundary $\mathcal{B}$, with the function $\gamma$, defined along $C_0^-$, the past light cone of $O$, the past end point of the central component $\mathcal{B}_0$ of $\mathcal{B}$, being unbounded, while the conclusion of Theorem 2.1 fails.

According to the definitions of [1] *initial conditions of bounded variation* means that

$$(4.1) \qquad \qquad \alpha = \frac{\partial}{\partial r}(r\phi) = \theta + \phi$$

is a function of bounded variation along $C_0^+$. This is equivalent to both $\theta$ and $\phi$ being functions of bounded variation along $C_0^+$. The total variation of $\phi$ along $C_0^+$ is equal to the integral

$$\int_0^\infty |\theta| \frac{dr}{r}$$

along $C_0^+$. In terms of the coordinates $(t, s)$ defined in Section 1 we have

$$(4.2) \qquad \qquad \text{T.V.}[\phi]_{t=0} = \int_{-\infty}^\infty |\theta(0, s)| ds.$$

Thus $\alpha_{t=0}$ being a function of bounded variation is equivalent to $\theta_{t=0}$ being a function of bounded variation which is integrable on the real line. Here, we shall characterize initial conditions in terms of the function $\theta_{t=0} = \vartheta$.

Suppose then that $\vartheta$ belongs to the set $\mathcal{E}$. Then according to Theorems 2.1 and 3.1, we have

$$(4.3) \qquad \qquad I(t) \to \vartheta(0) \ \text{ as } t \to \infty,$$

and with

$$h(s) = \sqrt{\frac{1}{s} \int_0^s (e^{s'} \vartheta(s') - \vartheta(0))^2 ds'}$$

$(s > 0)$,

$$(4.4) \qquad \limsup_{s \to 0+} \left\{ \frac{h(s)}{g(s)} \right\} < \infty,$$

where $g(s)$ is the function defined in Theorem 3.1. Let $f_1$ be a nonnegative integrable function on the real line, vanishing on $(-\infty, 0)$, whose restriction to $[0, \infty)$ is absolutely continuous and

$$(4.5) \qquad \lim_{s \to 0+} f_1(s) = 1.$$

Let also $f_2$ be a nonnegative integrable absolutely continuous function on the real line, vanishing on $(-\infty, 0]$, such that

$$(4.6) \qquad \limsup_{s \to 0+} \left\{ \frac{1}{g(s)} \sqrt{\frac{1}{s} \int_0^s e^{2s'} f_2^2(s') ds'} \right\} = \infty.$$

For example, we may define $f_2$ on $(0, 1)$ so that

$$\frac{1}{s} \int_0^s e^{2s'} f_2^2(s') ds' = g(s) \quad : \text{ for all } s \in (0, 1);$$

that is,

$$f_2(s) = e^{-s} \sqrt{\frac{d(sg(s))}{ds}} \quad : \text{ for all } s \in (0, 1).$$

Given real parameters $\lambda_1, \lambda_2$ we then consider the initial data given by:

$$(4.7) \qquad \tilde{\vartheta}_{(\lambda_1, \lambda_2)} = \vartheta + \lambda_1 f_1 + \lambda_2 f_2.$$

Since the restrictions of $\tilde{\vartheta}_{(\lambda_1, \lambda_2)}$ and $\vartheta$ to the interval $(-\infty, 0)$ coincide, the corresponding solutions coincide in the interior of $C_0^-$ (domain of dependence) and define the same functions $\gamma(t)$ and $I(t)$. Since by (4.3) and (4.5)

$$(4.8) \qquad \tilde{\vartheta}_{(\lambda_1, \lambda_2)}(0) = \lim_{t \to \infty} I(t) + \lambda_1$$

if $\lambda_1 \neq 0$, then Theorem 2.1 applies so that $\tilde{\vartheta}_{(\lambda_1, \lambda_2)} \notin \mathcal{E}$, while if $\lambda_1 = 0$, $\lambda_2 \neq 0$, then by (4.4) and (4.6),

$$\limsup_{s \to 0+} \left\{ \frac{\tilde{h}_{(\lambda_1, \lambda_2)}(s)}{g(s)} \right\} = \infty.$$

Hence Theorem 3.1 applies and again $\tilde{\vartheta}_{(\lambda_1, \lambda_2)} \notin \mathcal{E}$. Therefore the 2-dimensional linear subspace $\Pi_\vartheta$ of the space of initial data defined by

$$(4.9) \qquad \Pi_\vartheta = \{ \tilde{\vartheta}_{(\lambda_1, \lambda_2)} \ : \ (\lambda_1, \lambda_2) \in \Re^2 \}$$

intersects $\mathcal{E}$ at only one point, the point corresponding to $(\lambda_1, \lambda_2) = (0, 0)$, that is $\vartheta$ itself.

Suppose next that $\vartheta, \vartheta' \in \mathcal{E}$. Then

$$\Pi_\vartheta \bigcap \Pi_{\vartheta'} = \emptyset \quad \text{unless } \vartheta = \vartheta'.$$

For, if

$$(4.10) \qquad\qquad \tilde{\vartheta}_{(\lambda_1, \lambda_2)} = \tilde{\vartheta}'_{(\lambda_1', \lambda_2')}$$

then $\vartheta$ and $\vartheta'$ coincide on $(-\infty, 0)$; thus they define the same functions $\gamma(t), I(t)$ and $g(s)$, and the same functions $f_1(s), f_2(s)$. By (4.3),

$$\vartheta(0) = \vartheta'(0) = \lim_{t \to \infty} I(t).$$

Equality (4.10) at $s = 0$ then yields

$$\lambda_1 = \lambda_1'.$$

Thus (4.10) becomes

$$\vartheta - \vartheta' = (\lambda_2' - \lambda_2) f_2$$

and since by (4.4)

$$\limsup_{s \to 0+} \left\{ \frac{h(s)}{g(s)} \right\} < \infty, \quad \limsup_{s \to 0+} \left\{ \frac{h'(s)}{g(s)} \right\} < \infty$$

we obtain, in view of (4.6),

$$\lambda_2' = \lambda_2, \quad \vartheta = \vartheta'.$$

We have thus proved:

THEOREM 4.1.   *Consider the exceptional set $\mathcal{E}$ in the space of initial data* $BV \bigcap L^1$ *on the real line. Then for each $\vartheta \in \mathcal{E}$ there is a 2-dimensional linear subspace $\Pi_\vartheta$ such that*

$$(\Pi_\vartheta \setminus \{\vartheta\}) \bigcap \mathcal{E} = \emptyset.$$

*Moreover, if $\vartheta, \vartheta' \in \mathcal{E}$, then*

$$\Pi_\vartheta \bigcap \Pi_{\vartheta'} = \emptyset$$

*unless $\vartheta$ and $\vartheta'$ coincide. We may therefore say that $\mathcal{E}$ has positive codimension in the space of initial data.*

PRINCETON UNIVERSITY, PRINCETON, NJ
*E-mail address*: demetri@math.princeton.edu